\documentclass[reqno,12pt]{amsart}

\usepackage[margin=2.5cm]{geometry}
\usepackage[T1]{fontenc}
\usepackage[latin1]{inputenc}
\usepackage[english]{babel}
\usepackage{amsmath,amssymb,amsthm}
\usepackage{enumitem}
\usepackage{mymacros}
\usepackage{hyperref}
\usepackage{xcolor}
\usepackage{graphicx}

\DeclareMathOperator{\capa}{cap}

\renewcommand{\MR}[1]{}

\title{A note on simply interpolating sequences for the Dirichlet space}

\author{Nikolaos Chalmoukis}
\address{Fachrichtung Mathematik, Universit\"at des Saarlandes, 66123 Saarbr\"ucken, Germany}
\email{chalmoukis@math.uni-sb.de
}
\thanks{This research was supported by the Ministry of Science and Higher Education of the Russian Federation, agreement No. 075-15-2019-1619}

\DeclareMathOperator{\range}{range}
\begin{document}

\begin{abstract} 
We study simply interpolating sequences for the Dirichlet space in the unit disc. In particular we are interested in comparing three different sufficient conditions for simply interpolating sequences. The first one is the the so called one box condition, the second is the column bounded propererty for the associated Grammian matrix and the third one is a restricted version of the one box condition introduced by Bishop and, independently, by Marshall and Sundberg. We prove that the one box condition implies the column bounded property which in turn implies the restricted on box condition of Bishop-Marshall-Sundberg, and we give two counterexamples which show that the reverse implications fail even for weakly separated sequences.

\end{abstract}
\maketitle

\section{Introduction}

Let $\cH$ a {\it reproducing kernel Hilbert space} (rkHs) consisting of functions defined on some set $X$, with reproducing kernel $K:X\times X \to \bC $. Given a sequence of points $\cZ=\{z_n\} \subseteq X $ there exists a natural {\it weighted restriction operator} associated to the sequence $\{z_n\}$
\begin{align*}
    \cR_\cZ  :  \cH & \dashrightarrow \ell^2(\bN) \\
    f & \mapsto \{ K(z_n,z_n)^{-1/2} f(z_n) \}.
\end{align*} 
The dashed arrow indicates that in principle this operator does not necessarily take values in $\ell^2(\bN)$. The sequence $\cZ$ is called {\it universally interpolating} (UI) if $\range (\cR_\cZ) = \ell^2(\bN)$ and {\it simply interpolating} (SI) if $\range(\cR_\cZ) \supseteq \ell^2(\bN).$ Furthermore let $\cM(\cH)$ be the {\it multiplier algebra} associated to $\cH$ which consists of functions $m:X\to \bC$ such that $mf\in\cH$ whenever $f\in \cH.$ An interpolation problem can be also defined for for the multiplier algebra of $\cH$ by considering the (unweighted) restriction operator 
\begin{align*}
    \cT_\cZ: \cM(\cH) & \longrightarrow \ell^\infty(\bN) \\
    m & \mapsto \{ m(z_n) \}
\end{align*}
A sequence such that $\cT_\cZ(\cM(\cH))=\ell^\infty(\bN)$ it is called {\it multiplier interpolating} (MI). Such sequences have been first studied by Carleson for the multiplier algebra of the Hardy space $H^2(\bD)$, which coincides with the algebra of bounded analytic functions in $\bD$, in his seminal paper \cite{Carleson1958}. Shapiro and Shields have also studied universally and simply interpolating sequences for the Hardy space \cite{Shapiro61}. An interesting phenomenon which arises in this case is that all three types of interpolation describe the same class of sequences.

The fact that multiplier interpolating and universally interpolating sequences coincide turns out to be a common feature of many well known reproducing kernel Hilbert spaces. In particular it remains true in all spaces with the complete Nevanlinna Pick property \cite[Theorem 9.19]{AM00}. This is a technical property which we are not going to define here but the interested reader can consult the monograph \cite{AM00}. 
Many of the crucial properties of sequences regarding interpolation are encoded by the {\it Gram matrix}. Let $\cZ=\{z_n\}$ a sequence in $X$ and $\cH$ a rkHs of functions on $X$ with kernel $K$. The Gram matrix of $\cZ$ with respect to $K$ is defined as the infinite matrix 
\[ G=(g_{nm})_{n,m}:=\Bigg( \frac{K(z_n,z_m)}{ \sqrt{K(z_n,z_n)K(z_m,z_m)}} \Bigg)_{n,m}. \]

It is known that for a complete Nevanlinna Pick space universally interpolating (equivalently multiplier interpolating) sequences are characterized by the following two conditions on the Gram matrix 
\[ \sup_{n\neq m}|g_{nm}|<1, \quad \text{and} \quad G:\ell^2(\bN)\to\ell^2(\bN) \quad \text{is bounded}.\]
This is a highly non trivial theorem, proved first in \cite{Aleman17} as a consequence of the positive answer to the Kadison-Singer problem \cite{Marcus2015} and then in \cite{Hartz2020} using different methods. The first condition is usually called {\it weak separation} and the second {\it bounded Grammian condition.}

On the other hand the relation between simply and universally interpolating sequences is more subtle. In the best of our knowledge the only space in which it is known that there exist simply interpolating sequences which are not universally interpolating is the {\it Dirichlet space}.

The Dirichlet space $\cD$ is defined as the space of analytic functions $f$ in the unit disc such that 
\[ \int_{\bD}|f'(z)|^2dz \wedge d\bar{z} < + \infty. \]

It is well known that, once equipped with the inner product 
\[ \inner{f}{g}_\cD:=\frac{1}{\pi}\int_{\bD}f'(z)\overline{g'(z)}dz\wedge d\bar{z} + \frac{1}{2\pi}\sup_{0<r<1}\int_{\bT}f(r\theta)\overline{g(r\theta)}d\theta, \]
the space $\cD$ becomes a reproducing kernel Hilbert space with reproducing kernel given by the formula 
\[ K_\cD(z,w):=\frac{1}{z\bar{w}}\log\frac{1}{1-z\bar{w}}, \quad (z,w\in\bD). \]
Furthermore $\cD$ has the complete Nevanlinna Pick property and therefore multiplier and universally interpolating sequences coincide for $\cD$.  In two unpublished papers, \cite{Bishop94, Marshall94} Bishop and Marshall \& Sundberg gave a sufficient condition for a sequence to be simply interpolating for the Dirichlet space. In a subsequent work \cite[Corollary 5.1]{Be2002}, B\"oe gave another sufficient condition for simple interpolation which closely resembles the one of Bishop, Marshall \& Sundberg. The problem of simply interpolating sequences has been studied further in \cite{Arcozzi10} and \cite{ChalmAdv}.  

To take a closer look at these conditions let us introduce an auxiliary measure associated to a sequence of points $\cZ=\{z_n\}\subseteq \bD$. If $\delta_z$ is the Dirac measure at $z$, we define, 
\[ \mu_\cZ:=\sum_{n=1}^\infty\frac{\delta_{z_n}}{K_\cD(z_n,z_n)}. \]
Then, a sequence satisfies the bounded Grammian condition if and only if $\mu_\cZ$ is a {\it Carleson measure} for $\cD$. Carleson measures for the Dirchlet space have been characterized by Stegenga \cite{Stegenga80} as positive Borel measures in $\bD$ such that there exists a constant $C(\mu_\cZ)>0$ such that for all finite collections of arcs $I_1, I_2, \dots I_k $ in the unit circle, 
\[ \mu_\cZ\Big( \bigcup_{j=1}^k S(I_j) \Big) \leq C(\mu_\cZ) \capa \Big( \bigcup_{j=1}^k I_j \Big), \]
where, $S(I)$ is the {\it Carleson box} corresponding to $I$ defined by $S(I):=\{ z\in \bD\setminus \{0\} : z/|z|\in I, \, 1-|z| \leq |I| \}$ and $\capa$ is the classical {\it logarithmic capacity} in the complex plane. We present first the condition of  B\"oe which is easier to state. A sequence is simply interpolating if it is weakly separated and satisfies the {\it one box condition}, i.e., there exists a constant $C(\mu_{\cZ})>0$ such that 
\[ \mu_{\cZ}(S(I))\leq C(\mu_\cZ) \Big(\log\frac{1}{|I|} \Big)^{-1}, \quad (\forall I \subseteq \bT \,\, \text{arcs}). \tag{OB} \label{onebox}\]

Notice that the one box condition corresponds to Stegenga's condition for $k=1$. In particular if letting $I=\bT$ we find that the measure $\mu_\cZ$ must be finite. 

On the other hand the condition which was introduced by Bishop, Marshall \& Sundberg is a variant of the one box condition which we will call {\it restricted one box condition} and it involves testing on specific Carleson boxes. Let $\delta \in (0,1)$, for an arc $I\subseteq \bT$ we denote by $I^\delta$ the arc with the same midpoint as $I$ and length $|I|^\delta$. Let $\cZ=\{z_n\}$ be a sequence in the unit disc and $I_n$ the arc with midpoint $z_n/|z_n|$ and length $1-|z_n|$. The restricted one box condition requires that $\mu_\cZ$ is finite and there exist $C(\mu_\cZ)>0$ and $\delta\in(0,1)$ such that 
\[ \mu_\cZ(S(I_n^\delta)) \leq C(\mu_\cZ) \Big(\log\frac{1}{1-|z_n|} \Big)^{-1}, \quad (\forall n \geq 0). \tag{ROB} \]

It should be noted that both Bishop and Marshall-Sundberg construct  sequences which are not universally interpolating but are weakly separated and satisfy the restricted one box condition hence are simply  interpolating.

Yet another condition which appears in the literature  is what we call column bounded property of the Grammian. A sequence $\cZ$ is called {\it column bounded} (CB) if the corresponding Grammian $G$ is bounded as a linear operator from $\ell^2(\bN)$ to $\ell^\infty(\bN)$ or equivalently if the columns of $G$ form a bounded sequence in $\ell^2(\bN).$ It is the understanding of the author that so far the three conditions have fallen under the general term ``one box condition'' and at times have been considered to be equivalent. In this short note we will establish the exact implications between these three conditions.

\begin{thm} \label{thm:main}
Let $\cZ$ a sequence in the unit disc. The following implications hold
\[ \text{(OB)} \implies \text{(CB)} \implies \text{(ROB)}. \]
Furthermore, neither of the above implications is reversible even for sequences which are weakly separated.
\end{thm}

We should note at this point that there exist simply interpolating sequences which do not satisfy the restricted one box condition hence Theorem \ref{thm:main} does not include all simply interpolating sequences.

\subsection*{Notation}If $A$ and $B$ are quantities that depend on some parameters we will use the notation $A \lesssim B $ to denote that there exists $C>0$ such that $A\leq C B$ for all choices of the parameters. If $A\lesssim B$ and $B \lesssim A$ we write $A \approx B.$

\section{Proof of the main theorem}

Let us start by introducing some useful notation. For a point $z\in \bD\setminus\{ 0 \} $ we denote by $I(z)\subseteq \bT$ the unique arc such that $z^*:=z/|z|$ is the midpoint of $I(z)$ and $|I(z)|=1-|z|.$  
Furthermore let $\Gamma(z)$ be the Stolz angle centered at $z^*$, i.e. the set of points $w\in \bD$ such that $|z^*-w|\leq 2 (1-|w|)$.

Another object we are going to use are the dilated Carleson boxes. Let $\delta\in(0,1)$ a dilation factor, then we define the dilated Carleson box $S^\delta(I)$ as the Carleson box which corresponds to the interval with same midpoint as $I$ and length $|I|^\delta.$ Similarly we define $S^\delta(z):=S^\delta(I(z)).$

For proving the next lemma it will be convenient to have a more geometric description of the weak separation condition in the Dirichlet space. For this recall that the hyperbolic metric $d_\bD$ in the unit disc is given by the formula 
\[ d_\bD(z,w):= \frac{1}{2}\log\frac{1+|\varphi_z(w)|}{1-|\varphi_z(w)|}, \,\,\, \text{where} \,\,\, \varphi_z(w) = \frac{z-w}{1-\overline{z}w}. \] As it was shown in \cite{Bishop94} and \cite{Marshall94} a sequence $\cZ$ is weakly separated in the Dirichlet space if and only if there exists a constant $K=K(\cZ)$ such that 
\[ d_{\bD}(z_n,0)+1 \leq K d_{\bD}(z_n,z_m), \quad (\forall n\neq m). \]
 We will also denote by $B(z,r)$ the hyperbolic disc with center $z$ and radius $r$.
A convenient assumption we are going to make throughout the proof is that $|z_n|\geq 1/2, \forall z_n \in \cZ$.  This does not affect the statement of Theorem \ref{thm:main}  since the finite number of points in the disc $\{ |z|\leq 1/2 \}$ can be added or removed from the sequence without altering neither the assumption nor the conclusion of Theorem \ref{thm:main}.
With this at hand we can prove the following. 

\begin{lem}\label{lem1}
Let $\cZ \subseteq \bD$ be a sequence which satisfies the one box condition. Then for every $c>0$, there exists $N=N(c)\in \bN$ such that for every $z\in \bD$ the hyperbolic disc centered at $z$ of radius $ c(d_\bD(0,z)+1) $ contains at most $N$ points of the sequence $\cZ$.
\end{lem}

\begin{proof}
 The lemma follows by the following observation. Fix $c>0$, then there exists a positive constant $\delta=\delta(c)\in(0,1)$ such that the dilated Carleson box $S^\delta(z)$ contains the hyperbolic ball of radius $R=c(d_\bD(0,z)+1)$ centered at $z.$ We leave this elementary but tedious verification to the reader. 

Then, if $\mu_\cZ$ is the measure associated to the sequence
\[\mu_\cZ(B(z,R)) \leq \mu(S^\delta(z))  \lesssim \delta C(\mu_\cZ) \frac{1}{d_\bD(0,z)+1}. \]

One the other hand, 
\begin{align*} \mu_\cZ(B(z,r)) & = \sum_{z_n\in\cZ \cap B(z,R)}\frac{1}{K_\cD(z_n,z_n)} \approx \sum_{z_n\in \cZ \cap B(z,R)} \frac{1}{d_\bD(0,z_n)+1} \\ & \geq \frac{\#(\cZ\cap B(z,R))}{d_\bD(0,z)+R+1}  \gtrsim \frac{\#(\cZ\cap B(z,R))}{d_\bD(0,z)+1}.
\end{align*}

Hence there exists a natural number $N$ such that every hyperbolic disc $B(z,R)$ contains at most $N$ points.
\end{proof}

\begin{lem}
Let $\delta<1$ and $z\in \bD$. There exists $M=M(\delta) > 0$ such that the diameter the for any $k\in \bN$ and $z\in \bD\setminus\{0\}$ the set  $\Gamma(z)\cap(S^{\delta^{k+1}}(z) \setminus S^{\delta^k}(z))$ is contained in a hyperbolic ball of radius at most $M(1+d_\bD(0,z))$.
\end{lem}

\begin{proof}
Assume without loss of generality that $0<z<1$ and hence $\Gamma(z)=\Gamma(1)$. Let $E_k:= \Gamma(1)\cap(S^{\delta^{k+1}}(z) \setminus S^{\delta^k}(z))$. We have that
\begin{align*}d_\bD(0,1-(1-|z|)^{\delta^{k}})-d_\bD(0,1-(1-|z|)^{\delta^{k+1}}) \\
 = \delta^k\frac{1-\delta}{2}\log\frac{1}{1-|z|} + \cO(1)  \lesssim d_\bD(z,0)+1. \end{align*} 

For $ w=|w|e^{i\theta}\in \Gamma(1) $, by integrating the hyperbolic element of arc length over the curve $\gamma(s)=|w|e^{is}, 0\leq s \leq |\theta| $ we find
\[ d_\bD(w,|w|) \leq \int_\gamma \frac{ds}{1-s^2} = \frac{|\theta w|}{1-|w|^2} \lesssim \frac{|1-w|}{1-|w|} \leq 2 . \]
Finally if $w,\omega \in E_k$, 
\[ d_\bD(w,\omega) \leq d_\bD(w,|w|)+ d_\bD(|w|,|\omega|)+d_\bD(|\omega|,\omega) \lesssim d_\bD(z,0)+1.\] 
\end{proof}

\begin{proof}[Proof of Theorem \ref{thm:main}]
Suppose that $\cZ$ satisfies the one box condition condition and let $z_n \in \cZ$, we have the following elementary inequalities
\begin{align*}
    \sum_{m=1}^\infty |g_{nm}|^2 & = \sum_{m=1}^\infty \frac{
    \Big| \frac{1}{z_n\bar{z}_m}\log\frac{1}{1-z_n\overline{z}_m} \Big|^2
    }{K_\cD(z_n,z_n)K_\cD(z_m,z_m)} 
    \lesssim  \sum_{m=1}^\infty \frac{
    \Big( \log\frac{1}{|1-z_n\overline{z_m}|} \Big)^2
    }{K_\cD(z_n,z_n)K_\cD(z_m,z_m)} \\
   & =  \sum_{z_m\in S(z_n)}   \,\,\,\,\, + \sum_{z_m\not\in \Gamma(z_n)\cup S(z_n)} + \,\,\,\,\,\sum_{z_m\in\Gamma(z_n), z\not\in S(z_n)}   \frac{
    \Big( \log\frac{1}{|1-z_n\overline{z_m}|} \Big)^2 
    }{K_\cD(z_n,z_n)K_\cD(z_m,z_m)} \\
    & = (A) + (B)+(C), 
\end{align*}
where the quantity that we sum has the same occurrence. We proceed now to estimate the quantities $(A), (B)$ and $(C)$ separately. For $(A)$ we have

\begin{equation*}
    (A)  \approx \sum_{z_m\in S(z_n)}\frac{K_\cD(z_n,z_n)}{K_\cD(z_m,z_m)}=K_\cD(z_n,z_n)\mu_\cZ(S(z_n)) \lesssim \frac{K_\cD(z_n,z_n)}{\log\frac{1}{|I(z_n)|}} \lesssim 1,
\end{equation*}
where in the last inequality we have used that $\mu_\cZ$ satisfies the simple box condition of Bishop-Marshall-Sundberg. 
Let us now estimate (B),  for this fix $\delta = \frac{1}{2}$ and proceed as follows 

\begin{align*}
    (B) & \approx \sum_{z_m \not\in \Gamma(z_n)\cup  S(z_n)}  \dfrac{
    \Big( \log\frac{1}{|z_n^*-z_m^*|} \Big)^2
    }{K_{\cD}(z_n,z_n)K_{\cD}(z_m,z_m)} \\ 
    & = \frac{1}{K_{\cD}(z_n,z_n)}\sum_{k=0}^\infty \sum_{
    \overset{z_m \not\in \Gamma(z_n)\cup  S(z_n)}{z_m\in S^{2^{-(k+1)}}(z_n)\setminus S^{2^{-k}}(z_n)}
    } \Big( \log\frac{1}{|z_n^*-z_m^*|} \Big)^2 K_\cD(z_m,z_m)^{-1}  \\
   & \lesssim \frac{1}{K_\cD(z_n,z_n)} \sum_{k=0}^\infty \sum_{z_m\in S^{2^{-(k+1)}}(z_n)\setminus S^{2^{-k}}(z_n)} \Big( \log\frac{1}{(1-|z_n|)^{2^{-k}}} \Big)^2 K_\cD(z_m,z_m)^{-1} \\
   & \lesssim \frac{1}{K_\cD(z_n,z_n)} \sum_{k=0}^\infty \sum_{z_m\in S^{2^{-(k+1)}}(z_n)\setminus S^{2^{-k}}(z_n)} \Big[ K_\cD(z_m,z_m)^{-1} 2^{-2k} \Big( \log\frac{1}{1-|z_n|} \Big)^2 + K_\cD(z_m,z_m)^{-1} \Big] \\
   & \approx  K_\cD(z_n,z_n)  \sum_{k=0}^\infty 2^{-2k}\mu_\cZ(S^{2^{-(k+1)}}(z_n)) + \mu_\cZ(\bD) \\
   & \lesssim K_\cD(z_n,z_n) \sum_{k=0}^\infty 2^{1-k} \Big( \log\frac{1}{1-|z_n|} \Big)^{-1} + \mu_\cZ(\bD) 
    \lesssim 1.
\end{align*}
For estimating $(C)$ we apply Lemma \ref{lem1}, so we have

 \begin{align*}
     (C) & \approx \sum_{k=0}^\infty \sum_{z_m\in E_k} \frac{K_\cD(z_m,z_m)}{K_\cD(z_n,z_n)} \\
     & \leq \frac{N_0}{K_\cD(z_n,z_n)} \sum_{k=0}^\infty \log\frac{1}{(1-|z_n|)^{2^{-k}}}  \lesssim 1.
 \end{align*}

Suppose now that $\cZ$ is column bounded and $\eta \in (0,1)$. Note that

\begin{align*} \sum_{m=1}^\infty |g_{nm}|^2 & \gtrsim \sum_{z_m\in S^\eta(z_n)}  \frac{
    \Big| \log\frac{1}{(1-z_n\overline{z_m})} \Big|^2
    }{K_\cD(z_n,z_n)K_\cD(z_m,z_m)} \\
    & \approx \sum_{z_m\in S^\eta(z_n)}\frac{K_\cD(z_n,z_n)}{K_\cD(z_m,z_m)}=K_\cD(z_n,z_n) \mu_\cZ(S^\eta(z_n)).
    \end{align*} 

Notice that we can apply  \cite[Proposition 9.11]{AglerMcCarthy00} to conclude that we can write $\cZ$ as a finite union of weakly separated sequences\footnote{As it stated \cite[Proposition 9.11]{AglerMcCarthy00} applies to sequences such that the Grammian is bounded on $\ell^2(\bN)$ but its proof clearly only uses the $\ell^2(\bN) \to \ell^\infty(\bN)$ boundedness.} 

If a sequence is column bounded, then pick any $0<\delta<1$. By the column bounded property we have that 
\[ \sum_{z_m \in S^\delta(z_n)}  \frac{
    \Big( \log\frac{1}{|1-z_n\overline{z_m}|} \Big)^2 
    }{K_\cD(z_n,z_n)K_\cD(z_m,z_m)} \leq C.   \] 
    Then it just remains to notice that for $z_m \in S^\delta(z_n)$ we can estimate below  
    \[ \delta \log\frac{1}{1-|z_n|} \lesssim \log\frac{1}{|1-z_n\overline{z_m}|}. \]
  Where the implicit constant does not depend on $\delta$ or the sequence. Combining the two inequalities we get that 
  \[ \mu_\cZ(S^\delta(z_n)) \lesssim \frac{1}{\delta} \Big( \log\frac{1}{(1-|z_n|)^\delta} \Big)^{-1}.\]
  Which shows that the measure $\mu_\cZ$ is finite and the restricted one box condition holds for $\delta$.
\end{proof}

We shall construct now an example of a sequence which satisfies the column bounded property and is weakly separated but it does not satisfy the one box condition and an example of a weakly separated sequence which satisfies the restricted one box condition but not the column bounded property.

\begin{exa} 
Let $\frac{1}{2}<r<1, N \in \bN, 0< \ell < 1$. Consider furthermore an interval $I \subseteq \bT$ such that $ |I| = \ell$. On the arc $r I$ we place $N$ equispaced points and let $\cE = \cE(r,\ell, N)$ the set of these points. 

The one box condition is satisfied for this sequence with a constant which is at least 
\[ \Big( \log\frac{1}{|I|} \Big)\sum_{z\in \cE} \frac{1}{\norm{K_z}_\cD^2} \geq \frac{N \log\frac{1}{\ell}}{\log\frac{1}{1-r^2}}.  \]

On the other hand, if $z,w\in \cE$, 
\[ |K(z,w)| \approx \log \frac{1}{|1-z\overline{w}|} \approx \log \frac{1}{\max(1-|z|^2,1-|w|^2, |z^*-w^*|)} \leq \min \Big( \log\frac{1}{1-r^2}, \log \frac{N}{\ell} \Big).  \]

If we furthermore impose the conditions
\begin{equation}\label{assumptions} N^2 \leq \frac{\ell^2}{1-r^2}, \,\,\,\,\, N \leq \log\frac{1}{\ell} \end{equation}
we can conclude that 
\[ |K(z,w)| \lesssim \log \frac{1}{\ell}, \,\, z,w \in \cE. \]
Hence the $\ell^2$ norm of a column of the Grammian of the sequence $\cE$ can be estimated thusly
\[ \sum_{z\in \cE} \frac{|K(z,w)|^2}{ \norm{K_z}_\cD^2 \norm{K_w}_\cD^2} \lesssim N \Big( \frac{\log\frac{1}{\ell}}{\log\frac{1}{1-r^2}} \Big)^2. \]

Finally that under \eqref{assumptions} the sequence is weakly separated with e constant that does not depend on $r, N, \ell$ that is because the hyperbolic distance between two adjecent points in $\cE$ is comparable with absolute constants to $ \log \frac{\ell}{N(1-r^2)}$ which, by the first of the two equations in \eqref{assumptions}, is bigger than $\frac12\log\frac{1}{1-r^2}$.

To construct a sequence which is column bounded but it does not satisfy the one box condition consider $R = \log\frac{1}{1-r^2}$ as a free parameter and set 
\[ N = \lfloor (\log R )^2 \rfloor, \,\,\, \log\frac{1}{\ell} = \frac{R}{\log R}.   \]
It is then readily verified that the assumptions \eqref{assumptions} are satisfied for $R$ sufficiently big and 
\[ \frac{N \log\frac{1}{\ell}}{\log\frac{1}{1-r^2}} \approx \log R, \,\,\, N \Big( \frac{\log\frac{1}{\ell}}{\log\frac{1}{1-r^2}} \Big)^2 \approx 1. \]

To conclude the example it suffices to choose an increasing sequence of $R_n$ such that for the corresponding sequences of points $\cE_n$ that we have constructed before, such that the following holds 
\[  \frac{|K(z,w)|^2}{ \norm{K_z}_\cD^2 \norm{K_w}_\cD^2} < \frac{1}{ 2^{n} N_n \sum_{i=1}^{n-1} N_i }, z\in \cE_n, w\in \bigcup_{m<n} \cE_m. \] 
Then the union of the sequences $\cE_n$ is column bounded but it does not satisfy the one box condition.

To construct the second example, notice that as soon as $\cE$ is weakly separated then it satisfies the restricted one box condition with a constant that does not depend on $r,N,\ell$ since there exists some $0<\delta<1$ which depends only on the weak separation constant of $\cE$ such that $S^\delta(z)\cap S^\delta(w) = \empty$ for $z,w \in \cE$.

Consider now the following choice of parameters $R = \log\frac{1}{1-r^2}, N= 
\lfloor (\log R)^3 \rfloor, \log\frac{1}{\ell} = \frac{R}{\log R}$ and pick an increasing sequence of $R_n$ and choose the arcs $I_n$ in the construction of $\cE_n$ in such a way such that if $ z\in \cE_n, w\in \cE_m, $ then the dilated Carleson boxes $ S^\delta(z),  S^\delta (w)  $ are still disjoint.

Furthermore to satisfy the finite measure property it suffices to take $R_n$ such that 
\[ \sum_{n} \frac{(\log R_n )^3}{R_n } < \infty. \]
Then clearly the union of the sequences $\cE_n$ satisfies the restricted one box condition but not the column bounded property.

\end{exa}






\bibliographystyle{plain}
\bibliography{literature}
\end{document}